\documentclass[reqno]{article}
\usepackage[english]{babel}
\usepackage{enumerate,amsmath,amsthm,amsfonts,
amssymb,latexsym,mathrsfs,graphicx,graphics
}
\usepackage[all]{xy}
\usepackage[sort]{cite}

\usepackage[T1]{fontenc}
\usepackage{authblk}

\newtheorem{theorem}{Theorem}[section]
\newtheorem{proposition}[theorem]{Proposition}
\newtheorem{lemma}[theorem]{Lemma}
\newtheorem{corollary}[theorem]{Corollary}
\theoremstyle{definition}
\newtheorem{definition}[theorem]{Definition}
\newtheorem{exm}[theorem]{Example}

\newtheorem{prob}[theorem]{Problem}

\theoremstyle{remark}
\newtheorem{remark}[theorem]{Remark}

\newcommand{\ID}{\operatorname{\textsc{Id}}}
\newcommand{\vr}{^{\vee\odot}}
\renewcommand{\wr}{^{\wedge\oplus}}
\newcommand{\bydef}{\mathrel{\mathop:}=}
\newcommand{\ax}{\operatorname{Ax}}

\newcommand{\id}{\operatorname{id}}

\newcommand{\spec}{\operatorname{Spec}}
\newcommand{\spsr}{\operatorname{{}^{\sR\!}Spec}}
\newcommand{\spmv}{\operatorname{{}^{\MV\!}Spec}}
\newcommand{\Rad}{\operatorname{Rad}}
\newcommand{\Max}{\operatorname{Max}}

\newcommand{\luk}{\operatorname{\text{\L}}}

\newcommand{\fr}{\operatorname{Free}}
\newcommand{\Gab}{\cat{G}^\text{Ab}}

\newcommand{\supp}{\operatorname{supp}}

\newcommand{\gr}{Grothendieck }
\newcommand{\Id}{\operatorname{Id}}
\newcommand{\End}{\operatorname{End}}

\newcommand{\frab}{\operatorname{Free}_{\Gab}}
\newcommand{\sR}{{{}^s\!\mathcal{R}}}

\newcommand{\sL}{{{}^s\!\mathcal{L}}}

\newcommand{\MV}{\mathcal{MV}}

\newcommand{\Mod}{\textrm{-}{}^s\!\!\mathcal{M}\!\!\:\mathit{od}}

\newcommand{\cat}{\mathcal}

\newcommand{\B}{\operatorname{B}}

\newcommand{\R}{\mathbb{R}}
\newcommand{\N}{\mathbb{N}}

\renewcommand{\wp}{\mathscr{P}}

\renewcommand{\vec}{\mathbf}
\newcommand{\V}{\mathit{Var}}

\renewcommand{\phi}{\varphi}

\newcommand{\g}{\gamma}

\newcommand{\la}{\left\langle}
\newcommand{\ra}{\right\rangle}

\newcommand{\ov}{\overline}

\newcommand{\tensor}{\otimes}

\renewcommand{\P}{\operatorname{P_F}}

\def\amslatex\slash{{\protect\AmS-\protect\LaTeX}}

\begin{document} 

\title{{\bf The semiring-theoretic approach to MV-algebras: a survey}}

\author{Antonio Di Nola}
\affil{\small{\textit{Dipartimento di Matematica} \\ {\it Universit\`a di Salerno, Italy} \\ \texttt{adinola@unisa.it}}}

\author{Ciro Russo}

\affil{\small{\textit{Departamento de Matem\'atica} \\ {\it Instituto de Matem\'atica} \\ {\it Universidade Federal da Bahia, Salvador -- BA, Brazil} \\ \texttt{ciro.russo@ufba.br}}}

\maketitle
\date{}

\begin{abstract}
In this paper we review some of the main achievements of the semiring-theoretic approach to MV-algebras initiated and pursued mainly by the present authors and their collaborators. The survey focuses mainly on the connections between MV-algebras and other theories that such a semiring-based approach enabled, and on an application of such a framework to Digital Image Processing.

We also give some suggestions for further developments by stating several open problems and possible research lines.
\end{abstract}

\section{Introduction}
\label{intro}

\begin{quote}
A \emph{fuzzy set (class)} $A$ in $X$ is characterized by a \emph{membership (characteristic) function} $f_A(x)$ which associates with each point $x$ in $X$ a real number in the interval $[0,1]$, with value $f_A(x)$ at $x$ representing the ``grade of membership'' of $x$ in $A$.
\end{quote}

This was the original definition of fuzzy set given by Lotfi A. Zadeh \cite{zad} in 1965. The simplicity of such a definition is now well known to be only apparent; indeed many questions about fuzzy sets~--- especially those regarding their proper logical and algebraic settings~--- do not have ultimate answers and still nowadays are a matter of study.

Indeed, while classical set theory finds its natural algebraic framework in the theory of Boolean algebras and its proper logical setting in Classical Logic, the class of fuzzy subsets of a given set may have many different algebraic structures and corresponding logics. Concerning fuzzy logic, two main directions have to be distinguished (cf. \cite{zad2}): fuzzy logic in the broad sense and in the narrow one. The former serves mainly as a mathematical framework for fuzzy control, analysis of vagueness in natural language, and several other applications. It is developed more from a technical point of view, rather than a theoretical one, being in fact one of the techniques of soft-computing. Fuzzy logic in the narrow sense is a class of logical systems with a comparative notion of truth, whose theories are developed respecting the rules and the spirit of symbolic logic (syntax, semantics, axiomatization, truth-preserving deduction, completeness, etc.), thus being special cases of many-valued logics (cf. \cite{haj}).

In this paper we shall focus our attention on fuzzy logic in the narrow sense and in particular to the class of MV-algebras \cite{cha2}, the algebraic semantics of \L ukasiewicz propositional logic. Fuzzy logic in the narrow sense still has plenty of realizations in different formal systems; they are essentially the logics whose standard algebraic semantics is the real unit interval equipped with the usual order relation and a left-continuous t-norm which is meant to interpret the logical \emph{strong conjunction} (hence a \emph{strong intersection}, from a set-theoretic point of view) \cite{haj}.

\L ukasiewicz logic  \cite{luk1,luk2,luk3} is one of the longest-known many-valued logics. The algebraic semantics of such a logic, and connected structures, boasts a wide literature which is developed in (obviously term-equivalent) different fashions such as, for example, Wajsberg algebras, bounded hoops, and MV-algebras~--- the latter being the one we shall refer to (see, e.g., \cite{mv1,mv2,mv3,mv4,mv5,mv6,mv7,mv8,mv9,mv10,mv11,mv12,mv13,mv14,mun}). 

Among the various fuzzy logics (and corresponding algebras), \L ukasiewicz logic and MV-algebras are~--- in our opinion~--- the ones that best preserve some properties of symmetry. This statement can be better explained by the following observations.

In any Boolean algebra $\la B, \vee, \wedge, ^*, 0, 1 \ra$, the conjunction $\wedge$ is a residuated operation, thus it defines a residuum $\to$ which has the following properties: $x \to 0 = x^*$ and $x \to y = x^* \vee y$. Analogously, in any MV-algebra $\la A, \oplus, \odot, ^*, 0, 1 \ra$, the strong conjunction $\odot$ is residuated and its residuum $\to$ verifies similar equations: $x \to 0 = x^*$ and $x \to y = x^* \oplus y$.

The equation $x \to 0 = x^*$ also means that the logical equivalence between ``{\sc not} $\phi$'' and ``$\phi$ {\sc implies} {\sc False}'', which holds in classical logic and is intuitively well-founded and justified, is still verified in \L ukasiewicz fuzzy logic. Moreover, in \cite[page 28]{haj}, the author states:
\begin{quote}
Note that the dual notion of a t-conorm [\ldots] will not play any important role in this book. This is because conjunction and disjunction do not have any dual relation to the implication.
\end{quote}
In fact this is true for all of the best-known fuzzy logics except \L ukasiewicz one. Indeed, in \L ukasiewicz logic, (strong) conjunction and disjunction have, as standard algebraic semantics, a t-norm and its dual t-conorm~--- such a duality involving negation and implication too.

Last, De Morgan's laws for Boolean algebras have MV-algebraic analogues:
\begin{eqnarray*}
(x \vee y)^* = x^* \wedge y^* & \quad & (x \oplus y)^* = x^* \odot y^* \\ \nonumber
(x \wedge y)^* = x^* \vee y* & \quad & (x \odot y)^* = x^* \oplus y^*. \nonumber
\end{eqnarray*}

In \cite{dng}, the authors initiated a further approach to the study of MV-algebras by looking at these structures as a special class of idempotent semirings. This approach was eventually enforced in \cite{bdn} where the authors defined a category of idempotent semirings, called MV-semirings, proved that it is isomorphic to the one of MV-algebras, and furthermore characterized the class of commutative rings whose ideals form an MV-semiring.

The theories of semirings and idempotent semirings are nowadays well-established and find applications in many fields, such as discrete mathematics, computer science, computer languages, linguistic problems, finite automata, optimization problems, discrete event systems, computational problems (see, for instance, \cite{co0,gla,golan,gun,kat1,kat2,kat3,kat4,kol2,lit2,litmas}), and fuzzy sets too \cite{lit1}. The theory arising from the substitution of the field of real numbers with the ``max-plus'' (or ``tropical'') semifield $\la \R \cup \{-\infty\}, \max, +, -\infty, 0\ra$ is often referred to as \emph{idempotent} or \emph{tropical mathematics}. In this area, relevant works the reader may refer to are~--- among others~--- \cite{gun,rich}.

The aforementioned papers \cite{bdn,dng} made clear that looking at MV-algebras as special idempotent semirings gives the opportunity of importing results and techniques from the theory of semirings and, consequently, in some cases from ring theory too.

Therefore the main aim of subsequent studies such as, in particular, \cite{bdnf} and \cite{dnr2} was precisely to use the semiring-theoretic viewpoint as an inspiration and a tool for the study of MV-algebras.

The results achieved so far are very encouraging. Moreover, besides serving MV-algebra theory they suggest a possible ``payback,'' namely, that MV-algebras can  on their turn give ideas and tools to semiring and semifield theories. It is worth noticing also that, as well as MV-algebras, various other logic-related algebraic structures can be viewed as special idempotent semirings, and therefore this approach could be further extended.

In this paper we basically review the results achieved in \cite{bdnf}, \cite{dnr} and \cite{dnr2} by adopting the semiring-theoretic viewpoint on MV-algebras. Furthermore we pose several problems and questions with the aim of stimulating their investigation. The paper is organized as follows.

In Section \ref{mv-algebra} we recall the basic definitions and results on \L ukasiewicz logic and MV-algebras while Section \ref{semiring} contains preliminaries about semirings, MV-semirings, and their respective semimodules.

In Section \ref{group} we recall the characterization of finitely generated projective MV-semimodules that allows to construct of the \gr group $K_0 A$ of any MV-algebra $A$ and to prove that $K_0$ defines a functor from the category of MV-algebras to the one of Abelian groups.

Section \ref{semifield} is devoted to the connection between MV-semirings and idempotent semifields with a strong order unit ($u$-semifields, for short) which is basically a different presentation of Mundici's categorical equivalence between MV-algebras and Abelian $u\ell$-groups, the latter category being isomorphic to the one of idempotent $u$-semifields. Such a reformulation has two merits. On the one hand, MV-semirings and $u$-semifields are both subcategories of the one of idempotent semirings and, for any $u$-semifield $F$ and corresponding MV-semiring $A$, we have that the latter is also homomorphic image of the former in the larger category. On the other hand, such a quotient morphism from $F$ to $A$ induces a strong relationship between the respective categories of semimodules; more precisely, $A\Mod$ turns out to be, up to a natural isomorphism, a full subcategory of $F\Mod$.

In Section \ref{sheaf} we report the main construction and result of \cite{bdnf}, that is, the representation of any MV-semiring (and therefore of any MV-algebra) as the semiring of the global sections of its own \gr sheaf.

An application of MV-semirings and semimodules to Digital Image Processing was presented in \cite{dnr}. We briefly review that algorithm in Section \ref{luktran}.

Last, in Section \ref{conclusion} we add some concluding remarks and above all we pose some questions indicating directions for future research.

\section{\L ukasiewicz logic and MV-algebras}
\label{mv-algebra}

The language of \L ukasiewicz Propositional Logic (\L PL) consists of a denumerable set of variables $\V = \{x_i \mid i \in \omega\}$, the binary connective $\to$ and the unary one $\lnot$, and parentheses. Well-formed formulas are built recursively as usual:
\begin{enumerate}[(F1)]
\item every variable is a formula,
\item if $\phi$ is a formula, so is $\lnot\phi$,
\item if $\phi$ and $\psi$ are formulas, so is $\phi \to \psi$,
\item all formulas are built by iterative applications of (F1--F3).
\end{enumerate}
The only inference rule for \L PL is \emph{Modus Ponens}
$$\mathit{MP} \qquad \frac{\phi \qquad \phi \to \psi}{\psi},$$
and the set $\ax_{\luk}$ of axioms is determined by the following four axiom schemes:
\begin{enumerate}[(\L1)]
\item $\phi \to (\psi \to \phi)$,
\item $(\phi \to \psi) \to ((\psi \to \chi) \to (\phi \to \chi))$,
\item $((\phi \to \psi) \to \psi) \to ((\psi \to \phi) \to \phi)$,
\item $(\lnot \phi \to  \lnot \psi) \to (\psi \to \phi)$.
\end{enumerate}

An \emph{MV-algebra} is an algebra $\la A, \oplus, ^*, 0 \ra$ of type $(2,1,0)$ that satisfies the following equations
\begin{enumerate}[(MV1)]
\item $x \oplus (y \oplus z) = (x \oplus y) \oplus z$;
\item $x \oplus y = y \oplus x$;
\item $x \oplus 0 = x$;
\item $(x^*)^* = x$;
\item $x \oplus 0^* = 0^*$;
\item $(x^* \oplus y)^* \oplus y = (y^* \oplus x)^* \oplus x$.
\end{enumerate}

On every MV-algebra it is possible to define another constant $1$ and two further operations as follows:
\begin{itemize}
\item $1 = 0^*$,
\item $x \odot y = (x^* \oplus y^*)^*$,
\item $x \ominus y = x \odot y^*$,
\item $x \to y = x^* \oplus y$.
\end{itemize}
The following properties follow immediately from the definitions
\begin{enumerate}[(MV1)]
\setcounter{enumi}{6}
\item $1^* = 0$,
\item $x \oplus y = (x^* \odot y^*)^*$,
\item $x \oplus x^* = 1$,
\end{enumerate}
while (MV5) and (MV6) can be reformulated as follows
\begin{enumerate}[(MV1)]
\setcounter{enumi}{4}
\item $x \oplus 1 = 1$,
\item $(x \ominus y) \oplus y = (y \ominus x) \oplus x$.
\end{enumerate}

MV-algebras are known to be the equivalent algebraic semantics (in the sense of Blok and Pigozzi \cite{blpi}) of \L PL via the interpretation of each \L PL formula $\phi$ as the equation $\tau(\phi) \approx 0^*$ in the language of MV-algebras, where $\tau(\phi)$ is the MV-algebraic term in the same variables of $\phi$ obtained recursively by applying the following rules:
\begin{itemize}
\item $\tau(x_i) = x_i$, for any variable $x_i$,
\item $\tau(\psi \to \xi) = \tau(\psi)^* \oplus \tau(\xi)$ for any two formulas $\psi$ and $\xi$,
\item $\tau(\lnot \psi) = \tau(\psi)^*$ for any formula $\psi$.
\end{itemize}

It is well known that MV-algebras are naturally equipped with an order relation defined as follows
\begin{equation}\label{mvleq}
x \leq y \quad \textrm{ if and only if } \quad x^* \oplus y = 1.
\end{equation}
Moreover it is easy to verify that $x^* \oplus y = 1$ is indeed equivalent to each of the following conditions
\begin{itemize}
\item $x \odot y^* = 0$;
\item $y = x \oplus (y \ominus x)$;
\item there exists an element $z \in A$ such that $x \oplus z = y$.
\end{itemize}

The order relation also determines a lattice structure on $A$, with $0$ and $1$ respectively bottom and top element, and $\vee$ and $\wedge$ defined as follows
\begin{eqnarray}
&& x \vee y = (x \odot y^*) \oplus y = (x \ominus y) \oplus y, \label{mvvee} \\
&& x \wedge y = (x^* \vee y^*)^* = x \odot (x^* \oplus y). \label{mvwedge}
\end{eqnarray}
It is worth noticing that $\oplus$, $\odot$, and $\wedge$ distribute over any existing join and, analogously, $\oplus$, $\odot$ and $\vee$ distribute over any existing meet, in any MV-algebra $A$. In other words, for any family $\{b_i\}_{i \in I}$ of elements of $A$ for which there exists $\bigvee_{i \in I} b_i$, for any family $\{c_i\}_{i \in I}$ for which there exists $\bigwedge_{i \in I} c_i$, and for any $a \in A$, 
\begin{itemize}
\item[--] $a \bullet \bigvee_{i \in I} b_i = \bigvee_{i \in I} a \bullet b_i$, for $\bullet \in \{\oplus, \odot, \wedge\}$,
\item[--] $a \bullet \bigwedge_{i \in I} c_i = \bigwedge_{i \in I} a \bullet c_i$, for $\bullet \in \{\oplus, \odot, \vee\}$.
\end{itemize}

The best-known example of MV-algebra is the real unit interval $[0,1]$ with the sum $x \oplus y \bydef \min\{x + y, 1\}$ and the involution $x^* \bydef 1 - x$. The product is then defined by $x \odot y \bydef \max\{x + y - 1, 0\}$, and the resulting lattice structure is the natural totally ordered one. For a comprehensive exposition of MV-algebra theory, we refer the reader to \cite{mvbook}.

It is well known that the MV-algebra $[0,1]$ generates the whole class of MV-algebras both as a variety and as a quasi-variety. This means that, in order to determine whether an equation or a quasi-equation holds or not in all MV-algebras, it suffices to do it for $[0,1]$. Moreover, it is known as well that \L ukasiewicz propositional calculus is standard complete, that is, the following soundness and completeness theorem~--- which comes essentially from results by Rose and Rosser \cite{ros} and Chang \cite{cha,cha1,cha2}~--- holds.
\begin{theorem}\label{compthm}
Let $\phi$ be a well-formed formula of \L ukasiewicz propositional logic. The following conditions are equivalent:
\begin{enumerate}[(a)]
\item $\phi$ is provable in \textnormal{\L PL};
\item $\phi$ is valid in all MV-algebras;
\item $\phi$ is valid in all totally ordered MV-algebras;
\item $\phi$ is valid in the MV-algebra $[0,1]$.
\end{enumerate}
\end{theorem}

We also recall that, for any non-empty set $X$, $[0,1]^X$ with pointwise defined operations is an MV-algebra as well, and it is often referred to as the MV-algebra of fuzzy subsets of $X$ \cite{bel}.  

It is worth recalling the following representation theorem for MV-algebras \cite{dn1,dn2}.
\begin{theorem}\label{dinola}
Every MV-algebra $A$ is, up to an isomorphism, a subalgebra of $([0,1]^\star)^X$, the latter being the algebra of functions from some set $X$ to an ultrapower (which depends only on the cardinality of $A$) $[0,1]^\star$ of $[0,1]$.
\end{theorem}

We need also to recall some facts about ideals of MV-algebras. A subset $I$ of an MV-algebra $A$ is called an \emph{ideal} if it is a downward closed submonoid of $\la A, \oplus, 0\ra$, i.e. if it satisfies the following properties:
\begin{itemize}
\item $0 \in I$;
\item $I$ is downward closed, that is, for all $a \in I$ and $b \in A$, $b \leq a$ implies $b \in I$;
\item $a \oplus b \in I$ for all $a, b \in I$.
\end{itemize}
It is self-evident that $\{0\}$ and $A$ are ideals; an ideal $I$ is called \emph{proper} if $I \neq A$ or, that is the same, if $1 \notin I$. The set $\Id(A)$ of all ideals of an MV-algebra $A$ is partially ordered by set-inclusion and is closed under arbitrary intersections. For any subset $S$ of $A$, the \emph{ideal generated by $S$}, denoted by $(S]$, is defined as the intersection of all ideals of $A$ containing $S$; it is characterized by the following well-known result.
\begin{proposition}[cf. \cite{mvbook}]\label{idgen}
For any non-empty subset $S$ of $A$,
$$(S] = \{a \in A \mid a \leq x_1 \oplus \cdots \oplus x_n, \textrm{ for some } n \in \omega \textrm{ and } x_1, \ldots, x_n \in S\}.$$
\end{proposition}
An immediate consequence of Proposition \ref{idgen} is the following.
\begin{corollary}\label{properid}
Let $A$ be an MV-algebra and $S$ a subset of $A$. Then the ideal $(S]$ generated by $S$ is proper if and only if, for any $n \in \omega$ and for any $a_1, \ldots, a_n \in S$, $a_1 \oplus \cdots \oplus a_n < 1$.
\end{corollary}

It is also known that any non-trivial MV-algebra has maximal ideals; more precisely, any proper ideal of an MV-algebra is contained in a maximal one. The set of all maximal ideals of $A$ is denoted by $\Max A$, the intersection of all maximal ideals of $A$~--- which is, on its turn, an ideal~--- is called the \emph{radical} of $A$: $\Rad A \bydef \bigcap \Max A$. \emph{Semisimple} algebras, as usual, are defined as the subdirect products of simple algebras. However, in the theory of MV-algebras, they can be characterized as (non-trivial) algebras whose radical is $\{0\}$, and such a characterization is most often used as definition. It is worth noticing that Boolean algebras are all semisimple MV-algebras.

Let $P$ be an ideal of an MV-algebra $A$. If $a \wedge b \in P$ implies that $a \in P$ or $b \in P$, for any $a, b \in A$, then $P$ is called \emph{prime}. Usually the set of all prime ideals of $A$ is denoted by $\spec A$; here we shall denote it by $\spmv A$ in order to distinguish it from the set of prime ideals of the semiring reducts of $A$ which will be introduced in next section. $\spmv A$ can be naturally equipped with the Zarisky topology, namely with the topology generated by the basis of compact open sets $\{U(a) \mid a \in A\}$, where $U(a) \bydef \{P \in \spmv A \mid a \notin P\}$ for all $a \in A$. Such a topological space is called the \emph{spectral space} of the MV-algebra $A$.

Ideals and congruences of an MV-algebra $A$ are in one--one correspondence. Indeed, for any congruence $\sim$, $0/\sim$ is an ideal and, conversely, for any ideal $I$, the relation $\sim_I$ defined by ``$a \sim_I b$ iff $d(a,b) \bydef (a \odot b^*) \oplus (b \odot a^*) \in I$'' is a congruence on $A$~--- it is, in fact, the only one for which the class of zero is equal to $I$. Therefore, in MV-algebras, the congruence whose corresponding ideal is $I$ is often denoted by $I$ itself, and the congruence classes and the quotient algebra are denoted, respectively, by $a/I$ (for all $a \in A$) and $A/I$. The following result is a well-known characterization of congruence classes of quotient MV-algebras.
\begin{lemma}\label{classch}
Let $A$ be an MV-algebra and $I \in \Id(A)$. For all $a \in A$, $a/I = \{(a \oplus b) \odot c^* \mid b, c \in I\}$.
\end{lemma}

The bijective correspondence between ideals and congruences in MV-algebras implies also that \emph{simple} MV-algebras, i.e. those algebras whose congruence lattice is the two-element chain, have no non-trivial ideals, that is, $\Id(A) = \{\{0\}, A\}$. It may be worth recalling that ideals and filters are, in MV-algebras as well as in Boolean algebras, in one--one correspondence to each other. In this paper, following the tradition of MV-algebras, we shall deal with ideals, but all the results and constructions presented can be suitably reformulated in terms of filters. 

Besides the fact that Boolean algebras are MV-algebras, it is worth mentioning that the so-called \emph{Boolean elements} of any MV-algebra play an important role in the whole theory. Hereafter we recall their definition and basic properties.
\begin{itemize}
\item An element $a$ of an MV-algebra $A$ is called \emph{idempotent} or \emph{Boolean} if $a \oplus a = a$.
\item For any $a \in A$, $a \oplus a = a$ iff $a \odot a = a$.
\item An element $a$ is Boolean iff $a^*$ is Boolean.
\item If $a$ and $b$ are idempotent, then $a \oplus b$ and $a \odot b$ are idempotent as well; moreover we have $a \oplus b = a \vee b$, $a \odot b = a \wedge b$, $a \vee a^* = 1$ and $a \wedge a^* = 0$.
\item The set $\B(A) = \{a \in A \mid a \oplus a = a\}$ is a Boolean algebra, usually called the \emph{Boolean center} of the MV-algebra $A$.
\item For any $a \in A$ and $u \in \B(A)$, $a = (a \oplus u) \wedge (a \oplus u^*) = (a \odot u) \vee (a \odot u^*)$. 
\end{itemize}

\section{Idempotent semirings and semimodules}
\label{semiring}

In this section we recall some definitions and properties of idempotent semirings and semimodules over them, and we show how MV-algebras and semirings are connected. Besides the very basic notions (most of which can be found in \cite{golan}), we shall recall mainly material which will be needed in the next sections. Throughout the paper, by \emph{free} or \emph{projective} structure (semiring, semimodule, etc.) we shall always mean the corresponding categorical notion, i.e., respectively, a free or projective object in the corresponding category.

\begin{definition}\label{sr}
A \emph{(unital, idempotent) semiring} $\la S, \vee, \cdot, 0, 1 \ra$ is an algebraic structure with two binary operations and two constants such that
\begin{enumerate}[(S1)]
\item $\la S, \vee, 0\ra$ is a semilattice with identity,
\item $\la S, \cdot, 1\ra$ is a monoid,
\item $\cdot$ distributes over $\vee$ from both sides,
\item $a \cdot 0 = 0 = 0 \cdot a$ for all $a \in S$.
\end{enumerate}
A semiring $S$ is called \emph{commutative} if so is the multiplication, and a \emph{division semiring} if $\la S \setminus \{0\}, \cdot, 1\ra$ is actually a group. A commutative division semiring is called a \emph{semifield}.

Given two semirings $S$ and $S'$, a \emph{semiring homomorphism} from $S$ to $S'$ is any map $f: S \to S'$ that preserves the two binary operations and the two constants.
\end{definition}

Henceforth, by a ``semiring'' we shall always mean an idempotent unital semiring, except when differently specified.

\begin{definition}\label{sm}
Let $\la S, \vee, \cdot, 0, 1\ra$ be a semiring. An \emph{$S$-semimodule} is a bounded semilattice $\la M, \vee, 0 \ra$ with an external operation with coefficients in $S$, called \emph{scalar multiplication}, $\cdot: (a,x) \in S \times M \mapsto a x \in M$, such that the following conditions hold for all $a, b \in S$ and $x, y \in M$:
\begin{enumerate}[(SM1)]
\item $(a b) \cdot x = a \cdot (b \cdot x)$,
\item $a \cdot (x \vee y) = (a \cdot x) \vee (a \cdot y)$,
\item $(a \vee b) \cdot x = (a \cdot x) \vee (b \cdot x)$,
\item $0 \cdot x = 0 = a \cdot 0$,
\item $1 \cdot x = x$.
\end{enumerate}
\end{definition}
Let $S$ be a semiring and $M, N$ be two $S$-semimodules. A map $f: M \to N$ is an $S$-semimodule homomorphism if $f(x \vee y) = f(x) \vee f(y)$, for all $x, y \in M$, and $f(a x) = a f(x)$, for all $a \in S$ and $x \in M$.

If $S$ is a semiring and $X$ is a set, the \emph{support} of a map $f: X \to S$ is the set $\supp f = \{x \in X \mid f(x) \neq 0\}$. We have the following well-known result.
\begin{proposition}[cf. \cite{golan}]\label{freemod}
For any set $X$, the free $S$-semimodule $\fr_S(X)$ generated by $X$ is the set~--- denoted by $S^{(X)}$~--- of functions from $X$ to $S$ with finite support, equipped with pointwise join and scalar multiplication, and with the map $\chi: x \in X \mapsto \chi_x \in S^{(X)}$, where $\chi_x$ is defined, for all $x \in X$, by
\begin{equation}\label{chi}
\chi_x(y) = \left\{\begin{array}{ll} 0 & \textrm{if } y \neq x \\ 1 & \textrm{if } y = x\end{array}\right..
\end{equation}
In particular, for any natural number $n$, the $n$-generated free $S$-semimodule is up to isomorphisms $S^n$.
\end{proposition}
Obviously, every $S$-semimodule is homomorphic image of a free one.

In order to show next result, we introduce the following notation. Given a semiring $S$, let $S^{X \times (Y)}$ be the bounded semilattice of functions from $X \times Y$ to $S$~--- equipped with pointwise join~--- with finite support in the second variable, i.e. the functions $k: X \times Y \to S$ such that, for any fixed $x \in X$, the one-variable map $k(x,{}_\text{---}): Y \to S$ has finite support.

It is easy to see that $S^{X \times (Y)}$ also enjoys a structure of $S$-semimodule, in an obvious way. We present the next result along with its proof in order to better illustrate the situation. Moreover, the construction presented in the proof is useful also for the characterization of finitely generated projective semimodules.

\begin{theorem}\label{homofree}
Let $S$ be a semiring and $S^{(X)}$ and $S^{(Y)}$ free $S$-semimodules. The two $S$-semimodules $\hom_S(S^{(X)},S^{(Y)})$ and $S^{X \times (Y)}$ are isomorphic.
\end{theorem}
\begin{proof}
Let, for any $k \in S^{X \times (Y)}$,
\begin{equation}\label{transform}
\begin{array}{cccc}
h_k: & S^{(X)} & \to & S^{(Y)} \\
		 & f			 & \mapsto & \bigvee_{x \in X} f(x) k(x,{}_\text{---});
\end{array}
\end{equation}
$h_k$ is easily seen to be well-defined (since $f$ has finite support) and a semimodule homomorphism.

Conversely, let us observe that, for any $f \in S^{(X)}$, $f = \bigvee_{x \in X} f(x) \chi_x$, with the maps $\chi_x$ defined by (\ref{chi}). Then, for all $h \in \hom_S(S^{(X)},S^{(Y)})$,
$$h(f) = h\left(\bigvee_{x \in X} f(x) \chi_x\right) = \bigvee_{x \in X} f(x) h(\chi_x).$$
Hence $h = h_k$~--- as defined in (\ref{transform})~--- with $k: (x,y) \in X \times Y \mapsto h(\chi_x)(y) \in S$, and $k$ has, clearly, finite support in the second variable.

So let $\eta: k \in S^{X \times (Y)} \mapsto h_k \in \hom_S(S^{(X)},S^{(Y)})$; we shall now prove that $\eta$ is bijective. The fact that $\eta$ is surjective has just been proved. If $k \neq l \in S^{X \times (Y)}$, there exists a pair $(\ov x, \ov y) \in X \times Y$ such that $k(\ov x, \ov y) \neq l(\ov x, \ov y)$; then we have
$$h_k(\chi_{\ov x})(\ov y) = k(\ov x, \ov y) \neq l(\ov x, \ov y) = h_l(\chi_{\ov x})(\ov y),$$
whence $\eta$ is injective.

The fact that $\eta$ is a semimodule homomorphism is trivial.
\end{proof}

Now let us restrict our attention to endomorphisms of finitely generated free semimodules. For any semimodule $M$, the set $\End_S(M)$ of its $S$-semimodule endomorphisms has a natural structure of semiring with pointwise join, map composition, $0$-constant endomorphism and identity map. In what follows, by the \emph{endomorphism semiring} of a semimodule we shall mean the structure $\la \End_S(M), \vee, \cdot, 0, \id_M\ra$ where $\cdot$ is the composition in the reverse order, namely, $fg \bydef g \circ f$.

In order to present the representation theorem for finitely generated free semimodule endomorphisms, let us denote, for all $n \in \N$, by $M_{n}(S)$ the set of all $n \times n$ square matrices of elements of $S$. It is easy to verify that the structure $\la M_{n}(S), \vee, \star, o, \iota \ra$, where
\begin{itemize}
\item $o$ is the $0$-constant matrix,
\item $\iota$ is the matrix whose components are defined by $\id_{ij} = \left\{\begin{array}{ll} 1 & \textrm{if } i=j \\ 0 & \textrm{otherwise} \end{array}\right.$,
\item $\vee$ is the componentwise join,
\item the operation $\star$ is defined by $(a_{ij}) \star (b_{ij}) = \left(\bigvee_{k=1}^n a_{ik} b_{kj}\right)$,
\end{itemize}
is a semiring, called the \emph{semiring of $n \times n$ square matrices} over $S$.

\begin{theorem}[cf. \cite{dnr2}]\label{matrixendo}
The semirings $M_n(S)$ and $\End_S(S^n)$ are isomorphic, for any semiring $S$ and any natural number $n$.
\end{theorem}

\begin{theorem}\label{finproj}
An $n$-generated $S$-semimodule $M$ is projective if and only if there exist $\vec u_1, \ldots, \vec u_n \in S^n$ such that $M \cong S \cdot \{\vec u_i\}_{i=1}^n$ and the matrix $(u_{ij})$ is a multiplicatively idempotent element of the semiring $M_{n}(S)$.
\end{theorem}

Let $S = \la S, \vee, \cdot, 0,1 \ra$ be a semiring. An \emph{ideal} of $S$ is a nonempty subset $I \subseteq S$ such that
\begin{itemize}
\item $a \vee b \in I$ for all $a, b \in I$,
\item $ab \in I$ for all $a \in I$ and $b \in S$.
\end{itemize}
$I$ is proper if $I \neq S$ or, that is the same, if $1 \notin I$. An ideal $P$ of $S$ is \emph{prime} if $a \in P$ or $b \in P$ whenever $ab \in P$. As in the case of MV-algebras, the set $\spsr S$ of all prime ideals of a semiring $S$ is endowed with the Zarisky topology whose basis is $\{U(a) \mid a \in S\}$, where $U(a) \bydef \{P \in \spsr S \mid a \notin P\}$, for all $a \in S$.

Now we have all the necessary notions for linking MV-algebras and semirings. 
\begin{proposition}[cf. \cite{dng}]\label{mvsemi}
Let $A$ be an MV-algebra. Then $A\vr = \la A, \vee, \odot, 0, 1\ra$ and $A\wr = \la A, \wedge, \oplus, 1, 0 \ra$ are semirings. Moreover, the involution $^*: A \to A$ is an isomorphism between them.
\end{proposition}

Thanks to Proposition \ref{mvsemi}, we can limit our attention to one of the two \emph{semiring reducts of $A$}; therefore, whenever not differently specified, we will refer only to $A\vr$, all the results holding also for $A\wr$ up to the application of $^*$.

\begin{definition}[cf. \cite{bdn}]\label{luksemi}
An \emph{MV-semiring}, or \emph{\L ukasiewicz semiring}, is a commutative, additively idempotent semiring $\la A, \vee, \cdot, 0, 1 \ra$ for which there exists a map $^*: A \to A$~--- called the \emph{negation}~--- satisfying, for all $a, b \in A$, the following conditions:
\begin{enumerate}[(i)]
\item $ab = 0$ iff $b \leq a^*$ (where $\leq$ is naturally defined by means of $\vee$);
\item $a \vee b = (a^* \cdot (a^* \cdot b)^*)^*$.
\end{enumerate}
\end{definition}

\begin{proposition}[cf. \cite{bdn}]\label{luksemiprop}
If $A$ is an MV-semiring, we can define the operation $\oplus$ by
\begin{equation*}
a \oplus b = (a^* \cdot b^*)^* \quad \textrm{for all $a, b \in S$}.
\end{equation*}
Then the structure $A^\oplus = \la A, \oplus, ^*, 0 \ra$ is an MV-algebra.
\end{proposition}

Let $A$ be an MV-algebra. It is immediate to verify that $\spmv A$ do not need to coincide with $\spsr A\vr$. Moreover, the corresponding topological spaces need not be homeomorphic. We refer the reader to \cite{bdnf} for a detailed discussion about the relationship between such spaces.

By an \emph{MV-semimodule} we mean a semimodule on an MV-semiring. The presence of the involution $^*$ on MV-semirings allows the introduction of a special class of MV-semimodules. The defining property of such semimodules is basically a good behaviour of the scalar multiplication with respect to the MV-algebraic involution $^*$.

\begin{definition}\label{mvmod}
Let $A$ be an MV-semiring and $M$ an $A$-semimodule. $M$ is said to be a \emph{strong $A$-semimodule} provided it fulfils, for all $a, b \in A$, the following additional condition:
\begin{equation}\label{mvmodstar}
a \cdot x = b \cdot x \quad \textrm{for all } x \in M \qquad \textrm{implies} \qquad a^* \cdot x = b^* \cdot x \quad \textrm{for all } x \in M.
\end{equation}
\end{definition}

\begin{exm}\label{strong}
For any MV-algebra $A$, $\la A, \vee, 0 \ra$ is a strong $A\vr$-semimodule as well as $\la A, \wedge, 1  \ra$ is a strong $A\wr$-semimodule. It is easy to see also that any free MV-semimodule is strong.
\end{exm}
\begin{exm}\label{nonstrong}
Let $A$ be the MV-algebra $[0,1]$ and consider the join-semilattice $M = \la\left[0,\frac{1}{2}\right], \vee, 0 \ra$ as an $A$-semimodule with $\odot$ as the scalar multiplication. For any $a \leq 1/2$, $a \odot x = 0 \odot x$ for all $x \in M$ but, if we set for example $a = x = 1/2$, $a^* \odot x = 0$ while $0^* \odot x = 1/2$. Hence $M$ is not a strong MV-semimodule.
\end{exm}

\begin{proposition}[cf. \cite{dnr2}]\label{endmv}
Let $A$ be an MV-semiring and $M$ a semilattice with neutral element. Then $M$ is a strong $A$-semimodule if and only if $\End_{\sL}(M)$~--- which in general is not an MV-semiring~--- contains an MV-subsemiring that is homomorphic image of $A$ (in $\MV$).
\end{proposition}

The proof of the following proposition is trivial.
\begin{proposition}[cf. \cite{dnr2}]\label{}
The following hold for any MV-algebra $A$.
\begin{enumerate}[(i)]
\item For any MV-ideal $I$, the semilattice reduct of the quotient MV-algebra $A/I$ is a strong semimodule over $A$ with
$$\begin{array}{cccc}
\cdot: & A \times A/I & \to & A/I \\
			 & (a, x/I) & \mapsto & (a \odot x)/I.
\end{array}$$
\item If $B$ is an MV-algebra and $h \in \hom_{\MV}(A, B)$, every strong $B$-semimodule $N$ is a strong $A$-semimodule with
$$\begin{array}{cccc}
\cdot_A: & A \times N & \to & N \\
			 & (a, x) & \mapsto & h(a) \cdot_B x.
\end{array}$$
\end{enumerate}
\end{proposition}

We close this section by recalling the important relationship existing between lattice-ordered groups and idempotent division semirings which will be used in Section \ref{semifield}.

A \emph{division semiring} (neither necessarily idempotent nor commutative) is a semiring for which there exists a multiplicative inverse for each non-zero element. A commutative division semiring is called a \emph{semifield}.

A partially ordered group $\la G, \cdot, {}^{-1}, 1, \leq\ra$ is a group endowed with an order relation which is compatible with the binary operation, i.e., such that $a \leq b$ implies $ca \leq cb$ and $ac \leq bc$ for all $a,b,c \in G$. If the order relation defines a lattice structure, then the group is called a \emph{lattice-ordered group}, \emph{$\ell$-group} for short.

It is well known that a non-trivial $\ell$-group is necessarily infinite and unbounded; so let $\la G, \cdot, {}^{-1}, 1, \vee, \wedge\ra$ be an $\ell$-group and let us add a bottom element $\bot$ to $G$. If we set $x \cdot \bot = \bot = \bot \cdot x$ for all $x \in \ov G = G \cup \{\bot\}$, then the structure $\la \ov G, \vee, \cdot, {}^{-1}, \bot, 1 \ra$ is an idempotent division semiring. The same can be done by adding a top element and setting $\wedge$ instead of $\vee$ as the semiring sum; in this case we obtain the idempotent division semiring $\la G \cup \{\top\}, \wedge, \cdot, {}^{-1}, \top, 1 \ra$. 

Conversely, let $\la F, \vee, \cdot, {}^{-1}, \bot, 1\ra$ be an idempotent division semiring. So $\la F \setminus \{\bot\}, \cdot, {}^{-1}, 1\ra$ is a group and the semilattice order defined by $\vee$ is compatible with $\cdot$. Moreover, it is immediate to verify that, for all $x, y \in F \setminus \{\bot\}$, $x \wedge y = (x^{-1} \vee y^{-1})^{-1}$ and, therefore, $\la F \setminus \{\bot\}, \cdot, {}^{-1}, 1, \vee, \wedge\ra$ is a lattice-ordered group. 

We notice that the constructions above actually define a categorical isomorphism, and its inverse, between $\ell$-groups, with $\ell$-group homomorphisms, and idempotent division semirings with semiring homomorphisms.

\section{The \gr group: toward a $K$-theory of MV-algebras}
\label{group}

In the present section the construction of the $K_0$ group of unital semirings, with emphasis on MV-semirings, is shown, and it is proved that such a construction defines a functor from the category of unital semirings to the one of Abelian groups. As it is well known, algebraic $K$-theory is a deep and powerful part of homological algebra. Although semiring theory presents relevant differences with respect to ring theory, the construction of the $K_0$ functor is totally analogous to the classical case. Moreover, MV-algebras are, in some sense, more ``friendly'' semirings, and we believe that the study of $K_0$-invariants of MV-algebras may give new and unexpected information on the subject.

Let $S$ be a (not necessarily commutative) unital semiring, $\la\P(S), \oplus, [\{0\}]\ra$ the Abelian monoid of isomorphism classes of finitely generated projective (left) $S$-semimodules, and $J = \frab(\P(S))$ the free Abelian group generated by such isomorphism classes. For any finitely generated projective left $S$-semimodule $P$, we denote by $[P]$ its isomorphism class. Let $H$ be the subgroup of $J$ generated by all the expressions $[P] + [Q] - [P \oplus Q]$. The \emph{\gr group} $K_0S$ of a semiring $S$ is the factor group $J/H$.

Since the two semiring reducts of an MV-algebra $A$ are isomorphic and the \gr groups of isomorphic semirings are obviously isomorphic, it is possible to set the following:
\begin{definition}\label{k0}
The \gr group $K_0A$ of an MV-algebra $A$ is the one of either of its semiring reducts.
\end{definition}

\begin{lemma}[cf. \cite{dnr2}]\label{k0mfree}
For any semiring $S$, if we consider $\Gab$ as a concrete category over the one~--- $\cat M^{\text{Ab}}$~--- of Abelian monoids, $K_0S$ is $\cat M^{\text{Ab}}$-free over $\la \P(S), \oplus, [\{0\}]\ra$, with associated monoid morphism
\begin{equation}\label{k}
k_S: [P] \in \P(S) \mapsto [P]/H \in K_0S.
\end{equation}
\end{lemma}

The following result, which is the core of the main result of the section (Theorem \ref{k0thm}), is presented with its full proof in order to give the reader a better clue of how the $K_0$ functor works. 

\begin{lemma}[cf. \cite{dnr2}]\label{k0lemma}
Let $A$ and $B$ be two MV-algebras. Any MV-algebra homomorphism $f: A \to B$ induces a monoid homomorphism from $\P(A)$ to $\P(B)$.
\end{lemma}
\begin{proof}
By Theorem~\ref{finproj}, finitely generated projective semimodules over a semiring can be identified with multiplicatively idempotent square matrices with values in the same semiring. It is immediate to verify that, if $M \cong A \cdot (u_{ij})_{i,j=1}^m$ and $N \cong A \cdot (v_{ij})_{i,j=1}^n$ are finitely generated projective $A$-semimodules, the finitely generated projective $A$-semimodule $M \oplus N$ is isomorphic to $A \cdot (w_{ij})_{i,j=1}^{m+n}$ with
\begin{equation}\label{matrsum}
w_{ij} = \left\{\begin{array}{ll}
u_{ij} & \text{if } i,j \leq m \\
v_{(i-m)(j-m)} & \text{if } i,j > m \\
0 & \text{otherwise}
\end{array}\right..
\end{equation}
Moreover, if $M$ and $N$ are isomorphic, we can assume the corresponding matrices to have the same size. Indeed, suppose $m < n$, $M \cong M \oplus \underbrace{\{0\} \oplus \cdots \oplus \{0\}}_{n-m \text{ times}}$, hence the $m \times m$ matrix $(u_{ij})$ generates a semimodule isomorphic to the one generated by the $n \times n$ matrix $(u'_{ij})$ which coincides with $(u_{ij})$ on every entry $ij$ such that $i,j \leq m$ and is constantly equal to zero elsewhere.

Let now $(u_{ij})$ be an idempotent $n \times n$ $A$-matrix and consider the $B$-matrix $(f(u_{ij}))$. Since $f$ is an MV-homomorphism, it preserves all the MV-algebraic operations and the lattice structure, hence it is also a semiring homomorphism. So we have
$$\left(\bigvee_{k=1}^n f(u_{ik}) \odot f(u_{kj})\right) = \left(f\left(\bigvee_{k=1}^n u_{ik} \odot u_{kj}\right)\right) = (f(u_{ij})),$$
whence $(f(u_{ij})) \star (f(u_{ij})) = (f(u_{ij}))$ and $(f(u_{ij}))$ is an idempotent $n \times n$ $B$-matrix.

Now assume that $M \cong A \cdot (u_{ij})_{i,j=1}^n$ and $N \cong A \cdot (v_{ij})_{i,j=1}^n$ are isomorphic finitely generated projective $A$-semimodules. Then, for all $i = 1, \ldots, n$, there exist $a_{i1}, \ldots, a_{in}, b_{i1}, \ldots, b_{in} \in A$ such that $(u_{i1}, \ldots, u_{in}) = \bigvee_{k=1}^n a_{ik} \cdot (v_{k1}, \ldots, v_{kn})$ and $(v_{i1}, \ldots, v_{in}) = \bigvee_{k=1}^n b_{ik} \cdot (u_{k1}, \ldots, u_{kn})$. So each vector $f(\vec u_i) = (f(u_{i1}), \ldots, f(u_{in}))$ can be expressed as a linear combination of the vectors $\{f(\vec v_i)\}_{i=1}^n$ with the scalars $f(a_{i1}), \ldots, f(a_{in}) \in B$ and, conversely, each $f(\vec v_i)$ can be expressed as a linear combination of the vectors $\{f(\vec u_i)\}_{i=1}^n$ with the scalars $f(b_{i1}), \ldots, f(b_{in}) \in B$; this means that the subsemimodules of $B^n$ generated respectively by $\{f(\vec u_i)\}_{i=1}^n$ and $\{f(\vec v_i)\}_{i=1}^n$ are isomorphic.

The above guarantees that
\begin{equation}\label{fchap}
\begin{array}{cccc}
\hat f: & \P(A) & \to & \P(B) \\
		& \left[A \cdot (u_{ij})\right] & \mapsto & \left[B \cdot (f(u_{ij}))\right]
\end{array}
\end{equation}
is a well-defined map. The fact that $\hat f\left(\left[\{0\}\right]\right) = \left[\{0\}\right]$ is obvious. On the other hand, given two classes $\left[A \cdot (u_{ij})_{i,j=1}^m\right], \left[A \cdot (v_{ij})_{i,j=1}^n\right] \in \P(A)$, the semimodule $A \cdot (u_{ij}) \oplus A \cdot (v_{ij})$ is isomorphic to $A \cdot (w_{ij})_{i,j=1}^{m+n}$ with $(w_{ij})$ defined by (\ref{matrsum}), and
$$f(w_{ij}) = \left\{\begin{array}{ll}
f(u_{ij}) & \text{if } i,j \leq m \\
f(v_{(i-m)(j-m)}) & \text{if } i,j > m \\
f(0) = 0 & \text{otherwise}
\end{array}\right.,$$
whence
$$\hat f\left(\left[A \cdot (u_{ij}) \oplus A \cdot (v_{ij})\right]\right) = [B \cdot (f(u_{ij})) \oplus B \cdot (f(v_{ij}))] = \hat f\left(\left[A \cdot (u_{ij})\right]\right) \oplus \hat f\left(\left[A \cdot (v_{ij})\right]\right),$$
and the assertion is proved.
\end{proof}

Thanks to Lemma~\ref{k0lemma}, we can easily obtain the following result.
\begin{theorem}[cf. \cite{dnr2}]\label{k0thm}
$K_0$ is a functor from $\MV$ to $\Gab$.
\end{theorem}

\section{MV-algebras and idempotent semifields}
\label{semifield}

In this section we shall deal with the celebrated Mundici categorical equivalence between MV-algebras and lattice-ordered Abelian groups with a distinguished strong order unit. The main objective of the section is to show that looking at such an equivalence from the viewpoint of idempotent semiring theory gives new information and results. In order to better explain the latter statement, we need first to recall a few definitions which are necessary to our discussion. 

For a lattice-ordered Abelian group $G$, an element $u > 0$ is called a \emph{strong order unit} if, for all $x \in G$, $x > 0$, there exists $n \in \N$ such that $n u > x$. In the category $\ell\Gab_u$ of Abelian $\ell$-groups with a distinguished strong order unit the morphisms are $\ell$-group homomorphisms which send the distinguished strong unit of the domain to the one of the codomain.

In a well-known paper of 1986 \cite{mun}, Daniele Mundici proved that the category $\MV$ of MV-algebras, with MV-algebra homomorphisms, is equivalent to $\ell\Gab_u$, namely, the category of lattice-ordered Abelian groups with a distinguished strong order unit whose morphisms are lattice-ordered group homomorphisms that preserve the distinguished strong unit. The two functors that form such an equivalence are usually denoted by $\Gamma: \ell\Gab_u \to \MV$ and $\Xi: \MV \to \ell\Gab_u$; while the former is very easy to present and shall be recalled hereafter, the latter requires more work and the details of its construction are not really relevant to our discussion. However, a detailed yet relatively concise presentation of Mundici categorical equivalence is presented in \cite[Chapter 2]{mvbook}. 

Let $\la G, +, -, 0, \vee, \wedge, u \ra$ be an Abelian $u\ell$-group with distinguished strong order unit $u$. Then the MV-algebra $\Gamma(G)$ is $\la [0,u], \oplus, ^*, 0\ra$ with $x \oplus y \bydef (x + y) \wedge u$ and $x^* \bydef u - x$ for all $x, y \in [0,u]$. The mapping $\Gamma: G \in \ell\Gab_u \mapsto \Gamma(G) \in \MV$ is a full, faithful and isomorphism-dense functor.

As we observed in the previous sections, the category $\ell\Gab$ is isomorphic to the one of idempotent semifields; therefore, it is natural to call \emph{$u$-semifield} an idempotent semifield $\la F, \wedge, +, -, \top, 0, u \ra$ with an element $u >0$ ($u \neq \top$) such that, for all $x \in F$, $0 < x < \top$, there exists $n \in \N$ s.t. $n u > x$. So, given a $u$-semifield $\la F, \wedge, +, -, \top, 0, u \ra$, we obtain an MV-algebra by applying the $\Gamma$ functor to the Abelian $u\ell$-group $\la F\setminus\{\top\}, +, -, 0, \vee, \wedge, u\ra$, with $\vee$ defined by means of $\wedge$ and $-$. In what follows, given a $u$-semifield $\la F, \wedge, +, -, \top, 0, u \ra$, with a slight abuse of notation we shall denote by $\Gamma(F)$ the MV-algebra $\la [0,u], \oplus, ^*, 0\ra$. 

First of all, it must be noticed that we can now see Mundici's functors as an equivalence between two subcategories of the same category. Moreover, as we shall see, the functor $\Gamma$ naturally defines a semiring homomorphism from any $u$-semifield $F$ to its corresponding MV-semiring $\Gamma(F)$. Now, using well-known results on the tensor product of semimodules it is possible to show that this situation results in a full embedding of the category $\Gamma(F)\Mod$ into the category $F\Mod$. It is worth remarking also that the left-inverse of such an embedding somehow ``truncates'' $F$-semimodules to $\Gamma(F)$-semimodules similarly to how $\Gamma$ truncates idempotent $u$-semifields to MV-algebras. In order to present such results, we need to briefly recall some properties of tensor products of semimodules.

\begin{definition}\label{bimor}
Let $S$ be a semiring, $M$ and $N$ $S$-semimodule, and $L$ a bounded semilattice. A map $f: M \times N \to L$ is called an \emph{$S$-bimorphism} if, for all $x, x_1, x_2 \in M$, $y, y_1, y_2 \in N$, $a \in S$, the following conditions hold:
\begin{enumerate}[(i)]
\item $f(x_1 \vee x_2, y) = f(x_1,y) \vee f(x_2,y)$,
\item $f(x, y_1 \vee y_2) = f(x,y_1) \vee f(x,y_2)$,
\item $f(xa,y) = f(x,ay)$.
\end{enumerate}
The \emph{$S$-tensor product} $M \tensor_S N$ is the codomain of the universal bimorphism with domain $M \times N$.
\end{definition}

Now recall that for any set $X$ the free bounded semilattice over $X$ is the set $\wp_F(X)$ of the finite subsets of $X$ equipped with set-theoretic union and the bottom element $\varnothing$. For a semiring $S$, and $S$-semimodules $M$ and $N$, the tensor product $M \tensor_S N$ is, up to isomorphisms, the quotient $\wp_F(M \times N)/\equiv_R$ of the free semilattice generated by $M \times N$ with respect to the semilattice congruence generated by the set $R$:
\begin{equation}\label{R}
R = \left\{
	\begin{array}{l}
		\left(\left\{\left(\bigvee X, y\right)\right\}, \bigcup_{x \in X}\{(x,y)\}\right) \\
		\left(\left\{\left(x, \bigvee Y\right)\right\}, \bigcup_{y \in Y}\{(x,y)\}\right) \\
		\left(\{(x a, y)\}, \{(x,a y)\}\right) \\
	\end{array} \right\vert
	\left.
	\begin{array}{l}
	X \in \wp_F(M), y \in N \\
	Y \in \wp_F(N), x \in M \\
	a \in A \\
	\end{array}
	\right\}.
\end{equation}
The tensor product $M \tensor_S N$ naturally inherits a structure of $S$-semimodule from the ones defined on $M$ and $N$:
\begin{equation*}
\star: \ \left(a, \bigvee_{i=1}^n x_i \tensor y_i\right) \in S \times \left(M \tensor N\right) \ \to \ \bigvee_{i=1}^n (a \cdot x_i) \tensor y_i \in M \tensor N.
\end{equation*}

\begin{lemma}\label{indmod}
Let $S$ and $T$ be semirings and $h: S \to T$ a semiring homomorphism. Then $h$ induces a structure of $S$-semimodule on any $T$-semimodule. In particular, $h$ induces structures of $S$-semimodule $T$ itself.
\end{lemma}
\begin{proof}
Let $N$ be a $T$-semimodule with scalar multiplication $\cdot$. It is easy to verify that
\begin{equation}\label{starh}
\cdot_h: (a, x) \in S \times N \mapsto h(a) \cdot x \in N
\end{equation}
makes $N$ into an $S$-semimodule. Since $T$ is a semimodule over itself, the second part of the assertion follows immediately.
\end{proof}

The operation performed in (\ref{starh}) is well known in the theory of ring modules as \emph{restricting the scalars along $h$}. In fact it defines a functor
\begin{equation}\label{subh}
\begin{array}{cccc}
H: & B\Mod & \to & A\Mod \\
& N & \mapsto & N_h
\end{array}
\end{equation}
having both a right and a left adjoint. In particular, the left adjoint of $H$ is the functor
\begin{equation*}
\begin{array}{cccc}
H_l: & S\Mod & \to & T\Mod \\
			& M			& \mapsto & T \tensor_{S} M
\end{array}.
\end{equation*}
In the special case where $h$ is an onto semiring homomorphism we have
\begin{theorem}[cf. \cite{dnr2}]\label{emb}
Let $h: S \to T$ be an onto semiring homomorphism. Then the functor $H$ defined in (\ref{subh}) is a full embedding. Moreover the left adjoint $H_l$ is, up to a natural isomorphism, the left inverse of $H$, that is, $H_l \circ H$ and the identity functor $\ID_{T\Mod}$ are naturally isomorphic.
\end{theorem}

Now we can apply the constructions and results above to MV-algebras and $u$-semifields. As a first step, we observe that the functor $\Gamma$ defines a canonical onto semiring homomorphism from any idempotent $u$-semifield to its corresponding MV-algebra. 
\begin{lemma}[cf. \cite{dnr2}]\label{Gamma}
Let $F$ be an idempotent $u$-semifield. Then the function
\begin{equation*}
\begin{array}{cccc}
\gamma: & F & \to & \Gamma(F)\wr \\
			  & a 	& \mapsto & (a \vee 0) \wedge u
\end{array}
\end{equation*}
is a semiring onto homomorphism.
\end{lemma}

By Lemmas~\ref{indmod} and \ref{Gamma}, the homomorphism $\gamma$ defines an adjoint and coadjoint functor
\begin{equation}\label{G}
G: \Gamma(F)\Mod \to F\Mod
\end{equation}
for any idempotent $u$-semifield $F$. Combining Theorem \ref{emb} with Lemma \ref{Gamma} we obtain the following immediate result.
\begin{corollary}\label{mvemb}
The functor $G$ defined in (\ref{G}) is a full embedding and its left adjoint $G_l$ is its left inverse.
\end{corollary}

Now we shall explain what does it mean that the functor $G_l$ ``truncates'' $F$-semimodules to $\Gamma(F)$-semimodules similarly to how $\Gamma$ truncates idempotent $u$-semifields to MV-algebras. Let $\la F, \wedge, +, \top, 0, u \ra$ be an idempotent $u$-semifield, $A = \Gamma(F)$ and $F^{(X)}$ be the free $F$-semimodule over a given set $X$. Moreover, let us denote by $\chi_x$ and $\chi_x'$ the maps defined in (\ref{chi}) respectively for $F^{(X)}$ and $A^{(X)}$.

Let us consider the function
$$f: \ (a, \alpha) \in A \times F^{(X)} \ \mapsto \ a \oplus \bigwedge_{x \in \supp\alpha} \gamma(\alpha(x)) \oplus \chi_x' \in A^{(X)},$$
and let $\alpha, \alpha' \in F^{(X)}$ and $a,a'  \in A$. We have:
$$\begin{array}{l}
f(a \wedge a', \alpha) \\
= (a \wedge a') \oplus \left(\bigwedge\limits_{x \in \supp\alpha} \gamma(\alpha(x)) \oplus \chi_x'\right) \\
= \left(a \oplus \bigwedge\limits_{x \in \supp\alpha} \gamma(\alpha(x)) \oplus \chi_x'\right) \wedge \left(a' \oplus \bigwedge\limits_{x \in \supp\alpha} \gamma(\alpha(x)) \oplus \chi_x'\right) \\
= f(a, \alpha) \wedge f(a', \alpha),
\end{array}$$
similarly $f(a, \alpha \wedge \alpha') = f(a,\alpha) \wedge f(a,\alpha')$. Now let $b \in F$; if $b \neq \top$ then $\supp \alpha = \supp(b + \alpha)$ and we have
$$\begin{array}{l}
f(a, b + \alpha) \\
= a \oplus \bigwedge\limits_{x \in \supp\alpha} \gamma(b+\alpha(x)) \oplus \chi_x' \\
= a \oplus \bigwedge\limits_{x \in \supp\alpha} \gamma(b) \oplus \gamma(\alpha(x)) \oplus \chi_x' \\
= a \oplus \gamma(b) \oplus \bigwedge\limits_{x \in \supp\alpha} \gamma(\alpha(x)) \oplus \chi_x' \\
= f(a \oplus \gamma(b), \alpha).
\end{array}$$
If $b = \top$ then $f(a, \top + \alpha) = a \oplus 1 = 1 = f(a \oplus 1, \alpha) = f(a \oplus \g(\top), \alpha)$.

So $f$ is an $F$-bimorphism (see Definition \ref{bimor}), hence it defines a semilattice homomorphism $\phi: A \tensor_F F^{(X)} \to A^{(X)}$ which is actually an $A$-semimodule homomorphism for the commutativity of $A$. Let us now consider the map
$$\begin{array}{llll}
\psi: & A^{(X)} & \to 		 & A \tensor_F F^{(X)} \\
			& \alpha	& \mapsto & 0 \tensor \bigwedge\limits_{x \in \supp\alpha} \alpha(x) + \chi_x
\end{array}.$$
It is self-evident that $\phi \circ \psi = \id_{A^{(X)}}$; on the other hand, for any tensor $a \tensor \alpha \in A \tensor_F F^{(X)}$,
$$\begin{array}{l}
a \tensor \alpha = a \tensor \left(\bigwedge\limits_{x \in \supp\alpha} \alpha(x) + \chi_x\right) \\
= \bigwedge\limits_{x \in \supp\alpha} (a \tensor  (\alpha(x) + \chi_x)) \\
= \bigwedge\limits_{x \in \supp\alpha} ((a \cdot_\g \alpha(x)) \tensor \chi_x) \\
= \bigwedge\limits_{x \in \supp\alpha} ((a \oplus \g(\alpha(x))) \tensor \chi_x)
\end{array}$$
and
$$\begin{array}{l}
(\psi \circ \phi)(a \tensor \alpha) = (\psi \circ \phi)\left(\bigwedge\limits_{x \in \supp\alpha}(a \oplus \gamma(\alpha(x))) \tensor \chi_x\right) \\
= \bigwedge\limits_{x \in \supp\alpha}\left(a \oplus \gamma(\alpha(x)) \oplus (\psi \circ \phi)(0 \tensor \chi_x)\right) \\
= \bigwedge\limits_{x \in \supp\alpha} \left(a \oplus \gamma(\alpha(x)) \oplus \psi(0 \tensor \chi_x')\right) \\
= \bigwedge\limits_{x \in \supp\alpha}(a \oplus \gamma(\alpha(x))) \oplus \left(0 \tensor \chi_x\right) \\
= \bigwedge\limits_{x \in \supp\alpha} (a \oplus \gamma(\alpha(x))) \tensor \chi_x,
\end{array}$$
whence $\psi \circ \phi = \id_{A \tensor_F F^{(X)}}$. It follows that $A^{(X)}$ and $A \tensor_F F^{(X)}$ are isomorphic.

In the general case, if $M$ is an $F$-semimodule and $X$ is a set of generators for it, then $M$ is homomorphic image of $F^{(X)}$, that is, there exists an onto homomorphism $f: F^{(X)} \to M$. So, as in the previous case, we can define the map $f': (a, \alpha) \in A \times F^{(X)} \mapsto a \tensor f(\alpha) \in A \tensor_F M$ which is easily seen to be an onto $F$-bimorphism and, therefore, induces an onto $A$-semimodule homomorphism $\varphi': A \tensor_F F^{(X)} \to A \tensor_F M$. Hence $\varphi' \circ \psi$ is an $A$-semimodule onto homomorphism and $A \tensor_F M$ turns out to be homomorphic image of the free $A$-semimodule over the same set of generators $X$ via a sort of truncation of the original morphism $f: F^{(X)} \to M$.

\section{A semiring-based sheaf representation of MV-algebras}
\label{sheaf}

In the last decades, several sheaf representations of MV-algebras were presented (see, for example, \cite{dubpov,filgeo}). Usually, such representations use special classes of MV-algebras for the ``localization'' at each point of the underlying topological space. In this section we recall the sheaf representation of MV-semirings presented in \cite{bdnf}. Such a representation differs from all the previous ones exactly in the class of objects which locally represent the algebra at each point of the underlying space. Indeed, as we shall see, the main result presented in \cite{bdnf}, which is strongly based on Chermnykh's sheaf representation of semirings \cite{che}, represents each MV-algebra as a sheaf of idempotent semirings which, in general, may not be MV-semirings. 

Let $\la S, \vee, \cdot, 0, 1 \ra$ be a semiring and $D$ a submonoid of $\la S \setminus \{0\}, \cdot, 1\ra$. The semiring $S_D$ is constructed in the following way from $S$ and $D$.

Let $(a, b), (c, d) \in S \times D$ and define $(a, b) \sim (c, d)$ if and only if there exists $k \in D$ such that $adk = bck$. It is easy to verify that $\sim$ is an equivalence relation. Now let $S_D$ the quotient set $S \times D / \sim$, and let us denote by $a/b$ the equivalence class of the pair $(a, b)$, for all $a \in S$ and $b \in D$. By \cite[Proposition 42]{bdnf}, $S_D$ is a semiring with the following operations:
\begin{itemize}
\item $a/b \vee c/d \bydef (ad \vee bc)/bd$,
\item $(a/b) \cdot (c/d) \bydef ac/bd$,
\item the additive identity $0$ is the class $0/1$ and the unit $1$ is the class $1/1$.
\end{itemize}

Let $P \in \spsr S$ and set $D = S \setminus P$. $D$ is a multiplicative monoid and $S_D$ is a local semiring, that is, a semiring with a unique maximal ideal. In this situation we will write $S_P$ instead of $S_D$ in order to underline the role of the prime ideal $P$, and we will call such a semiring the \emph{localization of $S$ at $P$}.

\begin{remark}\label{notmv}
If $A$ is an MV-algebra and $P \in \spsr A\vr$, $(A\vr)_P$ need not be an MV-semiring, as shown in \cite[Example 45]{bdnf}.
\end{remark}

Let $S$ be a semiring, the \emph{\gr sheaf} of $S$ is the triple $G(S) = \left(\spsr S, E_S, \pi_S\right)$ where $E_S = \bigcup\{S_P \times \{P\} \mid P \in \spsr S\}$ and $\pi_S: E_S \to \spsr S$ is a local homeomorphism defined by $\pi_S(a/b, P) = P$. Henceforth we denote by $[s/t]_P$ the element $(s/t, P) \in E_S$ and by $\hat S$ the semiring of all \emph{global sections}, i.e. the semiring of the continuous functions of type $\hat s: P \in \spsr S \mapsto [s/1]_P \in E_S$, whose operations are defined by
\begin{itemize}
\item $\hat s \hat \vee \hat t \bydef \widehat{s \vee t}$, and
\item $\hat s \hat \cdot \hat t \bydef \widehat{s \cdot t}$,
\end{itemize}
with $\hat 0$ and $\hat 1$ as the respective identities.

The semiring-based sheaf representation of an MV-algebra relies on the following theorem due to Chermnykh.
\begin{theorem}[cf. \cite{che}]\label{cherm}
For any commutative semiring the map $\phi: s \in S \mapsto \hat s \in \hat S$ is a semiring isomorphism.
\end{theorem}

Then we have the following result.
\begin{theorem}[cf. {\cite[Theorem 50]{bdnf}}]
Any MV-algebra $A$ is isomorphic to the MV-algebra of all global sections of the \gr sheaf of its semiring reduct.
\end{theorem}

It is worth stressing that the \emph{stalks} of such a representation, i.e. the semirings $(A\vr)_P$, are not necessarily MV-semirings but only commutative idempotent semirings. Nonetheless, the algebra of all global sections is again an MV-semiring and therefore an MV-algebra.

\section{The \L TB algorithm for image compression}\label{luktran}

In the present section we present an application of MV-semimodules to Digital Image Processing. In \cite{dnr} the authors defined the \emph{fuzzy transform in \L ukasiewicz algebra} --- \emph{\L u\-ka\-sie\-wicz transform} for short --- as a particular homomorphism between semimodules over the MV-semiring $[0,1]$ and a special case of fuzzy transform \cite{ft} which was also implemented as an algorithm for compressing and reconstructing digital images. This research line is an extremely active one (see, e.g. \cite{hi1,hi2,hi3, pede, pefe, tso}), although the present authors did not continued working in it, and the \L ukasiewicz transform received the attention of several authors both from the area of Image Processing and the one of Fuzzy Equations (see, for example, \cite{dimloses, dimses, abb, ghko}). In order to make the presentation more appropriate for the present survey, we shall limit our attention to the case of finitely generated semimodules.  

According to Theorem \ref{homofree}, for any given $m,n \in \N$, a homomorphism between an $m$-generated and an $n$-generated free $[0,1]$-semimodule can be identified with an $m \times n$ matrix. Let $p \in [0,1]^{m \times n}$ be defined, for all $i = 1, \ldots, m$ and $j = 1, \ldots, n$, as follows
\begin{equation}\label{partition}
p(i,j)=
\begin{cases}
(n - 1)\frac{i-1}{m-1} - (j - 2) & \text{if $\frac{j-2}{n-1} \leq \frac{i-1}{m-1} \leq \frac{j-1}{n-1}$}\\
-(n - 1)\frac{i-1}{m-1} + j & \text{if $\frac{j-1}{n-1} \leq \frac{i-1}{m-1} \leq \frac{j}{n-1}$}\\
0 & \text{otherwise}
\end{cases}.
\end{equation}

Then the \L ukasiewicz transform $H_n$ from $[0,1]^m$ to $[0,1]^n$ is defined by 
\begin{equation}\label{hn}
H_{m,n}: \ (f_i)_{i+1}^m \in [0,1]^m \ \longmapsto \ \left(\bigvee_{i=1}^m f(i) \odot p(i,j)\right)_{j=1}^n \in [0,1]^n,
\end{equation}
and its inverse transform $\Lambda_{m,n}$ by
\begin{equation}\label{lambdan}
\Lambda_{m,n}: \ (g_j)_{j=1}^n \in [0,1]^n \ \longmapsto \ \left(\bigwedge_{j=1}^n p(i,j) \to g_j\right)_{i=1}^m \in [0,1]^m.
\end{equation}

It is worth noticing that, although the matrix $p$ was presented in its analytical form here (and restricted to the finite case), its definition actually comes from important syntactic objects of \L ukasiewicz logical calculus, as discussed in the original paper \cite{dnr}. Indeed, the one-variable functions $p_0, \ldots, p_{m-1} \in [0,1]^n$, defined by $p_i(j) = p(i+1,j)$, for $j=0, \ldots, m-1$, are the McNaughton functions corresponding to the critical separating class of formulas which allows to define normal forms in \L ukasiewicz logic. The interested reader may refer to \cite{dnlet}.

The application of the transform in order to compress a digital image proceeds, grossly speaking, in the following way.

Any (greyscale) image of $m \times n$ pixels can be formalised as a $[0,1]$-valued matrix, each entry of the matrix meaning the value (0 = black, 1 = white) of the pixel in the same position. Such a matrix is divided in submatrices of a fixed size $a \times b$, hence in subimages of $a \times b$ pixels ($a,b \in \N$). Each block is now a matrix $(a_{hk})_{(h,k) \in a \times b}$, and we rewrite such a matrix as a function from $m = ab$ to $[0,1]$, i.e., as an $m$-vector $\ov a$ as follows:
$$a_{hk} = \ov a_{(h-1)b+k}, \textrm{ for all $h = 1, \ldots, a$ and $k = 1, \ldots, b$}.$$
Then we choose two natural numbers $c < a$ and $d <b$ and set $n = cd$. Now we can apply the \L ukasiewicz transform $H_{m,n}$ to each $\ov a$ (and, possibly, rewrite $H_{m,n}(\ov a)$ back as a $c \times d$ matrix, namely, as an image of $c \times d$ pixels, thus composing a compressed image by rejoining the blocks). 

Then the reconstructed block will be $(\Lambda_{m,n} \circ H_{m,n})(\ov a)$, and the reconstructed image will be obtained by ``gluing'' together the reconstructed blocks. In the case of RGB colour images, we need to apply the same process to three matrices, each of which representing a ``channel,'' i.e., one of the three fundamental colours.

As observed in the original paper, there is a consistent gap between JPEG and \L TB in terms of performance. However, the latter possesses the interesting distinguishing feature of being lossless if applied to an image which has been previously compressed and reconstructed. In other words, the \L TB algorithm is lossy just for the first process and lossless for the subsequent ones. Moreover, \L TB showed to have a lower running time w.r.t. JPEG in most of the tested images. Besides that, JPEG is composed of two compression processes, one of which is lossless, while \L TB has a single, lossy, process.

\begin{figure}[tp]
\centering
	\includegraphics[width=2.30in]{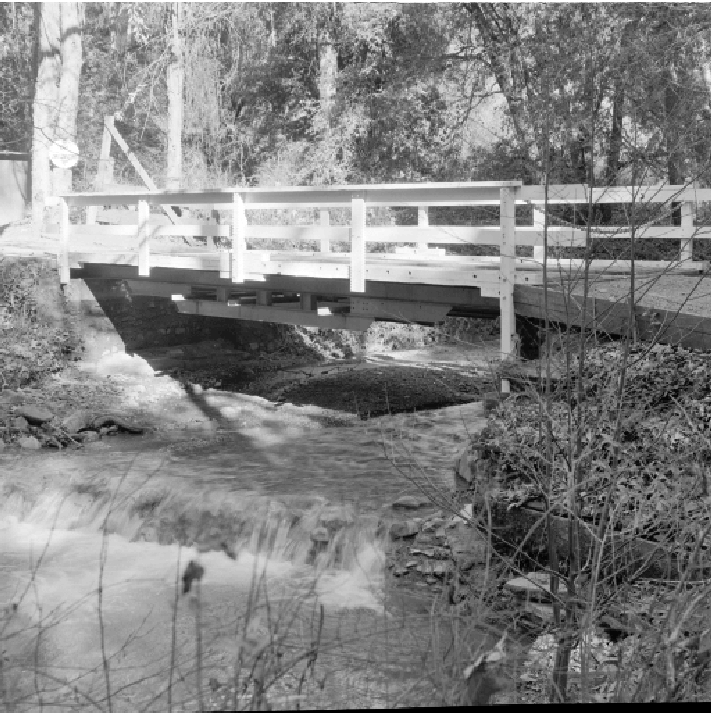}
		\caption{Bridge}
			\label{bridge}
\end{figure}

\begin{figure}[tp]
\centering
	\includegraphics[width=2.30in]{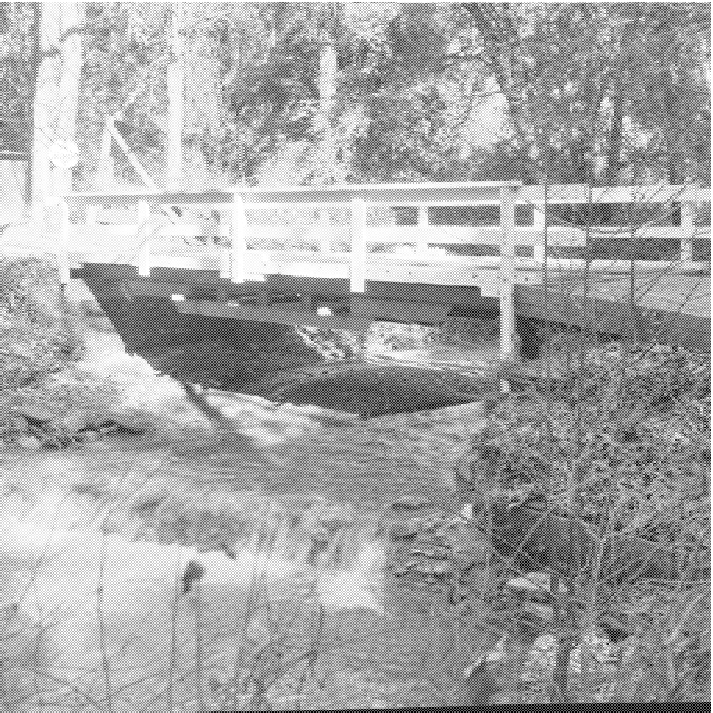}
	\caption{Bridge compressed and reconstructed by \L TB, $\rho = 0.5$}
			\label{22-21bridge}
\end{figure}

\begin{figure}[tp]
\centering
	\includegraphics[width=2.30in]{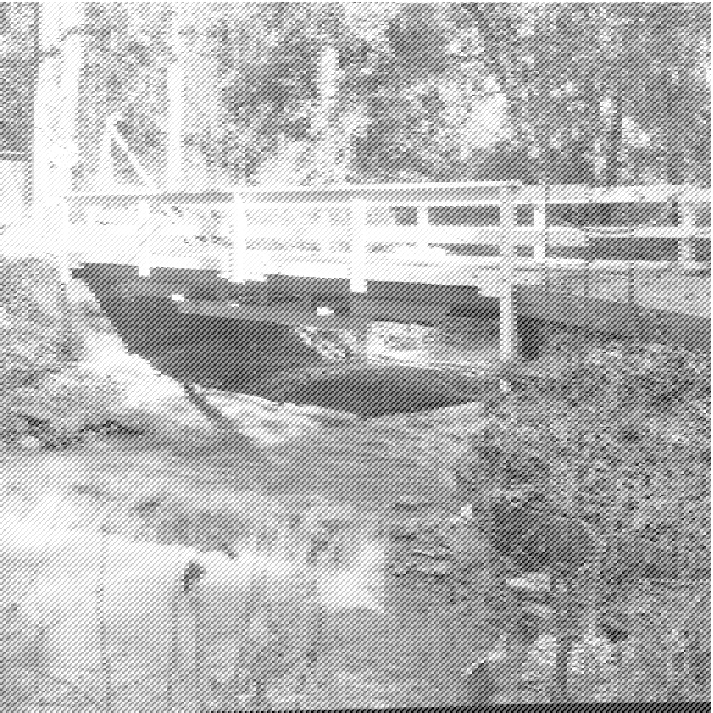}
		\caption{Bridge compressed and reconstructed by \L TB, $\rho = 0.25$}
			\label{44-22bridge}
\end{figure}

\begin{figure}[tp]
\centering
	\includegraphics[width=2.30in]{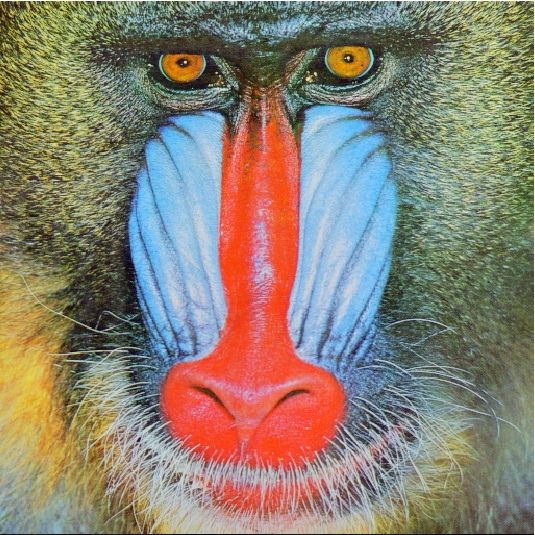}
\caption{Mandrill}
\label{mandrill}
\end{figure}

\begin{figure}[tp]
\centering
	\includegraphics[width=2.30in]{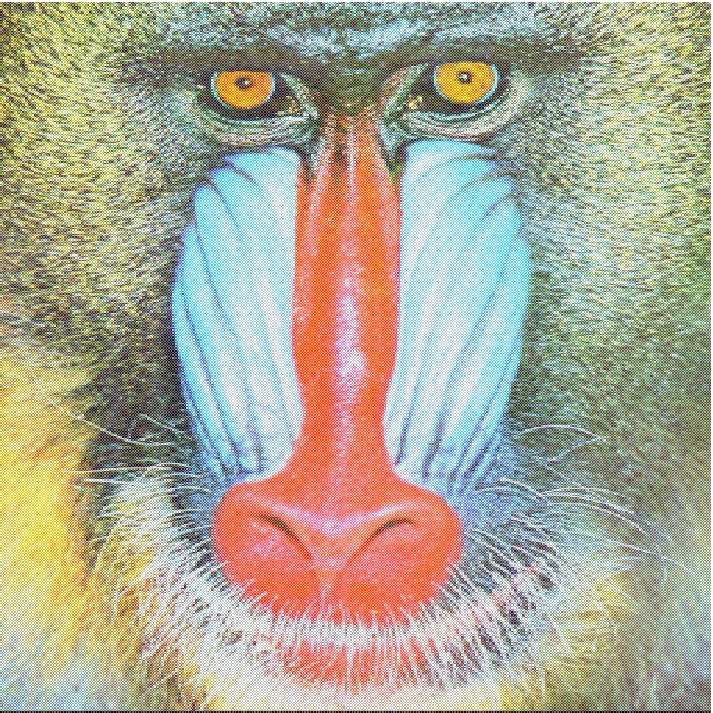}
\caption{Mandrill compressed and reconstructed by \L TB, $\rho = 0.5$}
\label{22-21mandrill}
\end{figure}

\begin{figure}[tp]
\centering
	\includegraphics[width=2.30in]{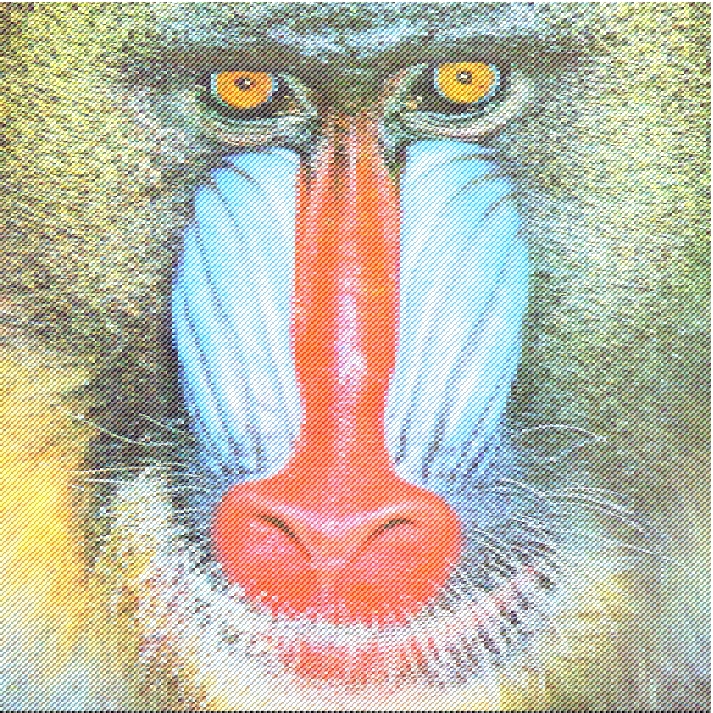}
\caption{Mandrill compressed and reconstructed by \L TB, $\rho = 0.25$}
\label{44-22mandrill}
\end{figure}

We show hereafter some concrete applications of the algorithm, the $\rho$ meaning the compression ratio
(see
Figs.~\ref{bridge}--\ref{44-22mandrill}).

\twocolumn

\begin{figure}[tp]
\centering
	\includegraphics[width=2.30in]{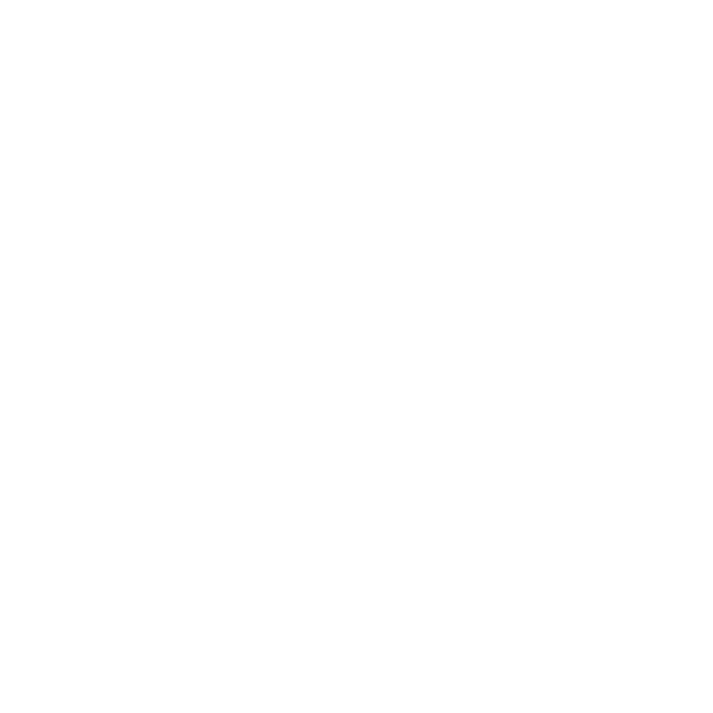}
\end{figure}

\twocolumn

\begin{figure}[tp]
\centering
	\includegraphics[width=2.30in]{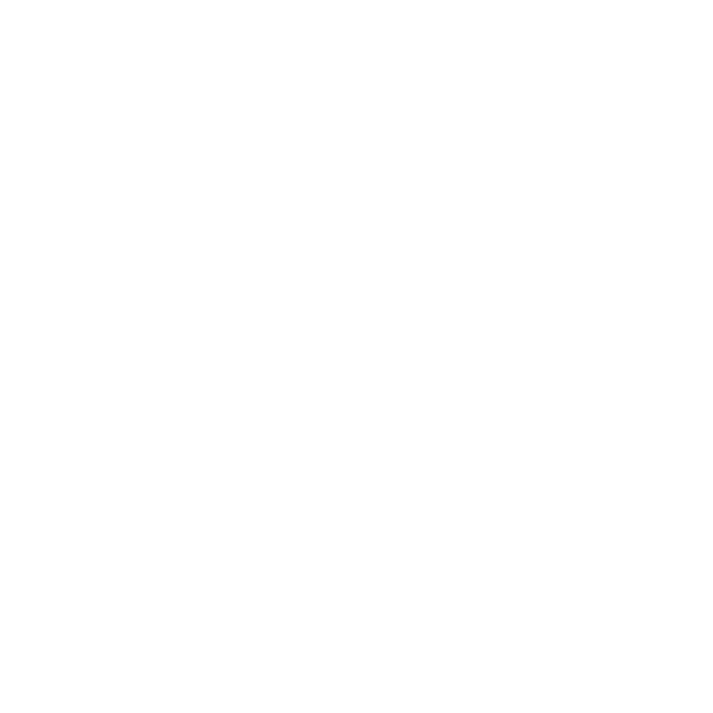}
\end{figure}

\onecolumn

\section{Conclusion}
\label{conclusion}

In this paper we reviewed some of the main results of the semiring-theoretic approach to the theory of MV-algebras initiated by Di Nola and B. Gerla \cite{dng} and developed by several researchers.

Even if the results overviewed in this paper are mostly of algebraic and categorical nature, several applications of MV-semirings have been already presented. The reader interested in this aspect may refer, for example, to \cite{dnr,rus0,rus,sch}. 

Moreover, we can say that the results achieved so far only scratched the surface on what can be done by applying semiring theory to MV-algebras. We conclude the paper with a series of open problems and suggestions for further investigations.

In Section \ref{group} we showed the construction of the \gr group of an MV-algebra. The functorial character of such a construction is easily proved by means of the matrix-based characterization of finitely generated projective semimodules (Theorem \ref{finproj}). Unfortunately, idempotent MV-algebraic matrices are almost completely unknown; therefore such a characterization of projective semimodules is not really easy to handle and, above all, does not really provide information on the $K_0$ group. Hence the following problems naturally arise.
\begin{prob}
Characterize idempotent matrices over MV-algebras.
\end{prob}
\begin{prob}
Describe the \gr groups of MV-algebras and find invariants under the $K_0$ functor.
\end{prob}

More in general, these problems seem to indicate also the need for a development of a linear algebra over MV-semirings.

The connection between MV-semimodules and semimodules over idempotent $u$-semifields established in Section \ref{semifield} let us catch a glimpse of a connection of the theory of MV-semimodules with the world of tropical geometry~--- the former being a sort of local version of the latter.
\begin{prob}
Which results and techniques, if any, can be borrowed from tropical geometry in order to obtain information on MV-algebras? Conversely, can MV-algebras make their own contribution to tropical geometry?
\end{prob}

As a whole, the results and constructions of Sections \ref{group} and \ref{semifield} show that the method of gathering information on a structure by studying its actions on external objects, as in classical ring theory, can be applied to MV-algebras if we look at them as special semirings. So, the  general question is obviously on how far this approach can lead us. In particular, the following problems can be immediately posed.
\begin{prob}
In \cite{kat4} the authors extended the notion of \emph{Morita equivalence} to semirings and presented several characterizations of Morita equivalent pairs of semirings. Basically, two semirings $S$ and $T$ are said to be Morita equivalent if the respective categories $S\Mod$ and $T\Mod$ of left semimodules are equivalent.

So it would be natural to look for a specific characterization of Morita equivalent pairs of MV-semirings. Moreover, the existence of a special subcategory of semimodules, namely, the one of strong semimodules, suggests a possible refinement of the notion of Morita equivalence and a specific investigation of such a concept.
\end{prob}

\begin{prob}
Develop a homological theory of MV-algebras starting from the general semiring-theoretic results of \cite{kat3}.
\end{prob}


\begin{thebibliography}{300}


\bibitem{abb}
Abbasi Molai A., Resolution of a system of the max-product fuzzy relation equations using $L\circ U$-factorization, {\em Informat. Sci.}, {\bf 234}, 86--96, 2013.

\bibitem{bel}
Belluce L. P., Semisimple algebras of infinite valued logic and bold fuzzy set theory, {\em Can. J. Math.}, {\bf 38}/6, 1356--1379, 1986.

\bibitem{bdn}
Belluce L. P., Di Nola A., Commutative Rings whose Ideals form an MV-algebra, {\em Mathematical Logic Quarterly}, {\bf 55}/5, 468--486, 2009.

\bibitem{bdnf}
Belluce L. P., Di Nola A., Ferraioli A. R., MV-semirings and their Sheaf Representations, {\em Order}, Published online Sept. 2011, doi:10.1007/s11083-011-9234-0.

\bibitem{blpi}
Blok W.~J., Pigozzi D., Algebraizable logics, {\em Memoirs of the Am. Math. Soc.}, {\bf 77}/396 (1989).

\bibitem{cha}
Chang C. C., Proof of an axiom of \L ukasiewicz, {\em Trans. Amer. Math. Soc.}, {\em 87}, 55--56, 1958.

\bibitem{cha1}
Chang C. C., Algebraic analysis of many valued logic, {\em Trans. Amer. Math. Soc.}, {\bf 88}, 467--490, 1958.

\bibitem{cha2}
Chang C. C., A new proof of the completeness of the \L ukasiewicz axioms, {\em Trans. Amer. Math. Soc.}, {\bf 93}, 74--90, 1959.


\bibitem{che}
Chermnykh V. V., Sheaf representations of semirings, {\em Russ. Math. Surv.}, {\bf 5}, 169--170, 1993.

\bibitem{mvbook}
Cignoli R. L. O., D'Ottaviano I. M. L., Mundici D., {\em Algebraic Foundations of Many-valued Reasoning}, Trends in Logic, {\bf 7}, Kluwer Academic Publishers, Dordrecht, 2000.

\bibitem{co0}
Cohen G., Gaubert S., Quadrat J.-P., Max-plus algebra and system theory: where we are and where to go now, {\em Annual Rev. Control}, {\bf 23}, 207--219, 1999.

\bibitem{dimloses}
Di Martino F., Loia V., Sessa S., Fuzzy transforms for compression and decompression of color videos, {\em Informat. Sci.}, {\bf 180} (20), 3914--3931, 2010.

\bibitem{dimses}
Di Martino F., Sessa S., Compression and decompression of images with discrete fuzzy transforms, {\em Informat. Sci.}, {\bf 177} (11), 2349--2362, 2007.

\bibitem{dn1}
Di Nola A., Representation and reticulation by quotients of MV-algebras, {\em Ricerche di Matematica}, {\bf 40}, 291--297, 1991.

\bibitem{dn2}
Di Nola A., MV-algebras in the treatment of uncertainty. In: L\"owen \& Roubens Eds., {\em Fuzzy Logic. Proceedings of the International Congress IFSA, Bruxelles 1991.}, 123--131, Kluwer, Dordrecht, 1993.

\bibitem{mv1}
Di Nola A., Dvure\v censkij A., State-morphism MV-algebras, {\em Ann. Pure Appl. Logic}, {\bf 161} (2), 161--173, 2009.

\bibitem{mv2}
Di Nola A., Flondor P., Leu\c stean I., MV-modules, {\em J. Algebra}, {\bf 267}(1), 21--40, 2003. 

\bibitem{dng}
Di Nola A., Gerla B., Algebras of \L ukasiewicz's Logic and their Semiring Reducts, {\em Idempotent Mathematics and Mathematical Physics}, 131--144, Contemp. Math., {\bf 377}, Amer. Math. Soc., 2005.

\bibitem{dnlet}
Di Nola A., Lettieri A., {\em On Normal Forms in \L ukasiewicz Logic},  Arch. Math. Logic, {\bf 43}, no. 6, 795--823, 2004.

\bibitem{dnr}
Di Nola A., Russo C., \L ukasiewicz Transform and its application to compression and reconstruction of digital images, {\em Informat. Sci.}, {\bf 177} (6), 1481--1498, 2007.

\bibitem{dnr2}
Di Nola A., Russo C., Semiring and semimodule issues in MV-algebras, {\em Comm. Alg.}, {\bf 41}/3 (2013), 1017--1048. \\ \texttt{arXiv:1002.4013v5 [math.RA]}

\bibitem{dubpov}
Dubuc E. J., Poveda Y.A., Representation theory of MV-algebras, {\em Ann. Pure Appl. Logic}, {\bf 161}(8), 1024--1046, 2010.

\bibitem{mv3}
Dvure\v censkij A., Effect algebras which can be covered by MV-algebras, {\em Internat. J. Theoret. Phys.}, {\bf 41}(2), 221--229, 2002.

\bibitem{mv4}
Dvure\v censkij A., Georgescu G., Iorgulescu A., Rudeanu S., Foreword: Multiple-valued logic and its algebras, {\em J. Mult.-Valued Logic Soft Comput.}, {\bf 16}(3--5), 219--220, 2010.

\bibitem{mv5}
Dvure\v censkij A., Kowalski T., Montagna F., State morphism MV-algebras, {\em Internat. J. Approx. Reason.}, {\bf 52}(8), 1215--1228, 2011.



\bibitem{filgeo}
Filipoiu A., Georgescu G., Compact and Pierce representations of MV-algebras, {\em Rev. Roumaine Math. Pures Appl.}, {\bf 40}(7--8), 599--618, 1995.

\bibitem{mv6}
Flaminio T., Montagna F., MV-algebras with internal states and probabilistic fuzzy logics, {\em Internat. J. Approx. Reason.}, {\bf 50}(1), 138--152, 2009.

\bibitem{mv7}
Georgescu G., Iorgulescu A., Pseudo-MV algebras, {\em G. C. Moisil memorial issue -- Mult.-Valued Log.}, {\bf 6}(1--2), 95--135, 2001.

\bibitem{ghko}
Ghodousian A., Khorram E., Linear optimization with an arbitrary fuzzy relational inequality, {\em Fuzzy Sets and Systems}, {\bf 206}, 89--102, 2012.

\bibitem{gla}
Glazek K., {\em A Guide to the Literature on Semirings and Their Applications in Mathematics and Information Sciences}, Kluwer Acad. Publ., Dordrecht, 2002.


\bibitem{golan}
Golan J. S., {\em Semirings and their Applications}, Kluwer Academic Publishers, Dordrecht, 1999.


\bibitem{gun}
Gunawardena J., An introduction to idempotency, In {\em Idempotency (Bristol, 1994)}, Vol. 11 of {\em Publ. Newton Inst.}, 1--49. Cambridge Univ. Press, 1998.

\bibitem{haj}
H\'ajek P., {\em Metamathematics of fuzzy logic}, Kluwer, 1998.

\bibitem{hi1}
Hirota K., Kawamoto K., Nobuhara H., Yoshida S.I., {\em On a lossy image compression/reconstruction method based on fuzzy relational equations}, Iran. J. Fuzzy Syst., {\bf 1}, 33--42, 2004.

\bibitem{hi2}
Hirota K., Nobuhara H., Pedrycz W., {\em Relational image compression: optimizations through the design of fuzzy coders and YUV color space}, Soft Computing, {\bf 9}, 471--479, 2005.

\bibitem{hi3}
Hirota K., Pedrycz W.,{\em Fuzzy relational compression}, IEEE Trans. Syst. Man Cyber.--Part B, {\bf 29} (3), 407--415, 1999.





\bibitem{kat1}
Katsov Y., Tensor products and injective envelopes of semimodules over additively regular semirings, {\em Algebra Colloq.}, {\bf 4}(2), 121--131, 1997.

\bibitem{kat2}
Katsov Y., On flat semimodules over semirings, {\em Algebra Univ.}, {\bf 51}, 287--299, 2004.

\bibitem{kat3}
Katsov Y., Toward homological characterization of semirings: Serre's conjecture and Bass's perfectness in a semiring context, {\em Algebra Univ.}, {\bf 52}, 197--214, 2004.

\bibitem{kat4}
Katsov Y., Nam T. G., Morita equivalence and homological characterization of semirings, {\em Journal of Algebra and its Applications}, {\bf 10}(3), 445--473, 2011.

\bibitem{kol2}
Kolokoltsov V. N., Maslov V. P., {\em Idempotent Analysis and Applications}, Kluwer Acad. Publ., Dordrecht, 1997.

\bibitem{mv8}
Kowalski T., Ono H., Fuzzy logics from substructural perspective, {\em Fuzzy Sets and Systems}, {\bf 161}(3), 301--310, 2010. 

\bibitem{mv9}
Leu\c stean I., Product MV-algebras and MV$^*$-algebras, {\em Information technology} (Bucharest, 1999), 1002--1005, Inforec, Bucharest, 1999.

\bibitem{lit1}
Litvinov G. L., Dequantization of Mathematics, idempotent semirings and fuzzy sets, In: E. P. Klement and E. Pap (eds.), {\em Mathematics of Fuzzy Systems, 25th Linz Seminar on Fuzzy Set Theory}, Linz, Austria, {\em Abstracts}, J. Kepler Univ., Linz, 2004.

\bibitem{lit2}
Litvinov G. L., Maslov dequantization, idempotent and tropical mathematics: a brief introduction, {\em Journal of Mathematical Sciences}, {\bf 140} (3), 2007.

\bibitem{litmas}
Litvinov G. L., Maslov V. P., Correspondence principle for idempotent calculus and some computer applications, Preprint IHES/M/95/33, Institut des Hautes \'Etudes Scientifiques, Bures-sur-Yvette 1995.

\bibitem{luk1}
\L ukasiewicz J., O zasadzie wy\l \c{a}czonego \'srodka. {\em Przegl\c{a}d Filozoficzny}, {\bf 13}, 372--373, 1910 (translation: On the Principle of the Excluded Middle. {\em History and Philosophy of Logic}, {\bf 8}, 67--69, 1987).

\bibitem{luk2}
\L ukasiewicz J., O logice tr\'ojwarto\'sciowej, {\em  Ruch Filozoficzny}, {\bf 5}, 170--171, 1920 (translation: On three-valued logic, in {\em Polish Logic, 1920--1939}, ed. S. McCall., 16--18).

\bibitem{luk3}
\L ukasiewicz J., Philosophische Bemerkungen zu mehrwertigen Systemen des Aussagenkalk\"uls, {\em Comptes rendus de la Société des Sciences et des Lettres de Varsovie}, cl. III, {\bf 23}, 51--77, 1930 (translation: Philosophical remarks on many-valued systems of propositional logic, in {\em Polish Logic, 1920--1939}, ed. S. McCall., 40--65).


\bibitem{mv10}
Moisil G. C., Sur les logiques de \L ukasiewicz \`a un nombre fini de valeurs, {\em Rev. Roumaine Math. Pures Appl.}, {\bf 9}, 905--920, 1964.

\bibitem{mv11}
Moisil G. C., Sur le calcul des propositions de type sup\'erieur dans la logique de Lukasiewicz \`a plusieurs valeurs, {\em An. Univ. Bucure?ti Ser. Acta Logica}, {\bf 10}, 93--104, 1967.

\bibitem{mv12}
Montagna F., Panti G., Adding structure to MV-algebras, {\em J. Pure Appl. Algebra}, {\bf 164} (3), 365--387, 2001.


\bibitem{mv13}
Montagna F., Subreducts of MV-algebras with product and product residuation, {\em Algebra Universalis}, {\bf 53} (1), 109--137, 2005.


\bibitem{mun}
Mundici D., Interpretation of $\mathnormal{AF \ C}^*$-algebras in \L ukasiewicz sentential calculus, {\em J. Functional Analysis}, {\bf 65}, 15--63, 1986.

\bibitem{mv14}
Mundici D., {\em Advanced \L ukasiewicz calculus and MV-algebras}, Trends in Logic, {\bf 35}, Springer, 2011.

\bibitem{ft}
Perfilieva I., Fuzzy Transforms: Theory and Applications, {\em Fuzzy Sets and Systems}, {\bf 157}(8), 993--1023, 2006.

\bibitem{pede}
Perfilieva I., De Baets B., Fuzzy transforms of monotone functions with application to image compression, {\em Informat. Sci.}, {\bf 180} (17), 3304--3315, 2010.

\bibitem{pefe}
Petrosino A., Ferone A., Rough fuzzy set-based image compression, {\em Fuzzy Sets and Systems}, {\bf 160}(10), 1485--1506, 2009.

\bibitem{rich}
Richter-Gebert J., Sturmfels B., Theobald T., First Steps in Tropical Geometry, In: G. L. Litvinov and V. P. Maslov (eds.), {\em Idempotent Mathematics and Mathematical Physics}, Contemp. Math., {\bf 377}, 289--318, Amer. Math. Soc., Providence, Rhode Island, 2005.

\bibitem{ros}
Rose A., Rosser J. B., Fragments of many-valued sentential calculi, {\em Trans. Amer. Math. Soc.}, {\bf 87}, 74--80, 1958.

\bibitem{rus0}
Russo C., A Unified Algebraic Framework for Fuzzy Image Compression and Mathematical Morphology, In: V. Di Gesù, S.K. Pal, and A. Petrosino Eds.: WILF 2009, {\em LNAI}, {\bf 5571}, 205--212, 2009. \\ \texttt{arXiv:1002.0982 [cs.IT]}

\bibitem{rus}
Russo C., Quantale Modules and their Operators, with Applications, {\em Journal of Logic and Computation}, {\bf 20}/4, 917--946, Oxford University Press, 2010. \\ \texttt{arXiv:1002.0968v1 [math.LO]}

\bibitem{sch}
Schwarz S., \L ukasiewicz Logic and Weighted Logics over MV-Semirings, {\em J. Autom. Lang. Comb.}, {\bf 12} (4), 485--499, 2007.

\bibitem{tso}
Tsolakis D., Tsekouras G. E., Niros A. D., Rigos A., On the systematic development of fast fuzzy vector quantization for grayscale image compression, {\em Neural Networks}, {\bf 36}, 83--96, 2012.

\bibitem{zad}
Zadeh L. A., Fuzzy sets, {\em Inform. and Control}, {\bf 8}, 338--353, 1965.

\bibitem{zad2}
Zadeh L. A., Preface, in: R. J. Marks II (ed.), {\em Fuzzy logic technology and applications}, IEEE Publications, 1994.


\end{thebibliography}
\end{document}